\documentclass[reqno,centertags]{amsart}
\usepackage{eurosym}
\usepackage{amsmath}
\usepackage{amssymb}
\usepackage{amsfonts}

\newtheorem{theorem}{Theorem}
\theoremstyle{plain}

\newtheorem{lemma}{Lemma}

\newtheorem{problem}{Problem}

\numberwithin{equation}{section}

\input{tcilatex}

\begin{document}
\title[]{Tubes of finite $II$-type Gauss map}
\author{Hassan Al-Zoubi }
\address{Department of Mathematics, Al-Zaytoonah University of Jordan, P.O.
Box 130, Amman, Jordan 11733}
\email{dr.hassanz@zuj.edu.jo}
%\author{Shadi Al-Zu'bi}
%\address{Al-Zaytoonah University of Jordan, Department of Computer Science }
%\email{smalzubi@zuj.edu.jo}
%\author{Stylianos Stamatakis}
%\address{Department of Mathematics, Aristotle University of Thessaloniki}
%\email{stamata@math.auth.gr}

%\author{Hani Almimi}
%\address{Al-Zaytoonah University of Jordan, Department of Computer Science}
%\email{Hani.Mimi@zuj.edu.jo}
\date{}

\begin{abstract}
In this paper, we continue the classification of finite type Gauss map surfaces in the 3-dimensional Euclidean space $\mathbb{E}^{3}$. We present an important family of surfaces, namely, tubes in $\mathbb{E}^{3}$. We show that the Gauss map of a tube is of infinite type corresponding to the second fundamental form.
\end{abstract}

\maketitle

\section{Introduction}

The theory of surfaces of finite Chen type was introduced by B.-Y. Chen in the late 1970s. From that moment on, interest in this field became widespread by many differential geometers, and it remains until this moment.

Let $M^{2}$ be a surface in the Euclidean space $\mathbb{E}^{3}$. Denote by $\Delta^{I}$ the Beltrami-Laplace operator on $M^{2}$ with respect to the first fundamental form $I$ of $M^{2}$. Then $M^{2}$ is said to be of finite type, if the position vector $\boldsymbol{x}$ of $M^{2}$ can be written as a finite sum of nonconstant eigenvectors of the operator $\Delta^{I}$, that is if
\begin{equation*}
\boldsymbol{x}=\boldsymbol{x}_{0}+\sum_{i=1}^{k}\boldsymbol{x}_{i},\ \ \ \ \
\ \Delta ^{I}\boldsymbol{x}_{i}=\lambda _{i}\boldsymbol{x}_{i},\ \ \
i=1,...,k,  \label{0}
\end{equation*}
where $\boldsymbol{x}_{0}$ is a fixed vector and $\boldsymbol{x}_{1},...,\boldsymbol{x}_{k}$ are nonconstant maps such that $\Delta ^{I}\boldsymbol{\
x}_{i}=\lambda _{i}\boldsymbol{x}_{i},i=1,...,k$.

In \cite{C5} B.-Y. Chen mentions the following problem

\begin{problem}
\label{(1)}Determine all surfaces of finite Chen $I$-type in $\mathbb{E}^{3}$.
\end{problem}

Regarding the above problem, important families of surfaces were studied by different authors by proving that finite type ruled surfaces \cite{C7}, finite type quadrics \cite{C8}, finite type tubes \cite{C4}, finite type cyclides of Dupin \cite{D1} and finite type spiral surfaces \cite{B1} are surfaces of the only known examples in $\mathbb{E}^{3}$. However, for other classical families of surfaces, such as surfaces of revolution, translation surfaces as well as helicoidal surfaces, the classification of its finite type surfaces is not known yet.

In the framework of surfaces of finite type S. Stamatakis and H. Al-Zoubi in \cite{S1} introduced the notion of surfaces of finite type with respect to the second or third fundamental form. For results concerning surfaces of finite type with respect to the second or third fundamental form (see \cite{A3}, \cite{A4}, \cite{A5}, \cite{A6}, \cite{A7}, \cite{A9}, \cite{A11}).

Another generalization can be made by studying surfaces in $\mathbb{E}^{3}$ of coordinate finite type, that is, their position vector $\boldsymbol{x}$ satisfying the relation
\begin{equation}
\Delta ^{J}\boldsymbol{x}=A\boldsymbol{x}, \ \ \ J = II, III,  \label{5}
\end{equation}%
where $A\in \mathbb{Re}^{3\times 3}$.

In \cite{A1} H. Al-Zoubi and S. Stamatakis studied ruled and quadric surfaces of coordinate finite type with respect to the third fundamental form satisfying
\begin{equation}
\Delta ^{III}\boldsymbol{x}=A\boldsymbol{x}.   \label{6}
\end{equation}

The same authors in \cite{S2} classified the class of surfaces of revolution satisfying condition (\ref{6}). Later in \cite{A2} the class of translation surfaces was studied, and it was shown that Sherk's surface is the only translation surface satisfying (\ref{6}).
B. Senoussi and M. Bekkar, in \cite{S0} classified the Helicoidal surfaces with $\triangle^{J}r = Ar, J = I, II, III$.

On the other hand, H. Al-Zoubi, T. Hamadneh in \cite{A14} classified the class of surfaces of revolution satisfying
\begin{equation}
\Delta ^{II}\boldsymbol{x}=A\boldsymbol{x}.   \label{7}
\end{equation}
where $A\in \mathbb{Re}^{3\times 3}$.

Similarly, we can apply the notion of surfaces of finite type to a smooth map, for example, the Gauss map of a surface. Here in this kind of research, many results can be found in (\cite{A8}, \cite{A10}, \cite{A12}, \cite{A13}, \cite{B2}, \cite{B3}, \cite{B4}, \cite{B5}, \cite{B6}, \cite{K1}).

In this paper, we will pay attention to surfaces of finite $II$-type. Firstly, we will give a formula for $\Delta^{II}\boldsymbol{x}$. Further, we continue our study by proving finite type Gauss map for an important class of surfaces, namely, tubes in $\mathbb{E}^{3}$.

\section{Preliminaries}

Let $\boldsymbol{x} = \boldsymbol{x}(u^{1}, u^{2})$ be a regular parametric representation of a surface $M^{2}$ in the Euclidean 3-space $\mathbb{E}^{3}$ referred to any coordinate system, which does not contain parabolic points. We denote by $b_{ij}$ the components of the second fundamental form $II=b_{ij}du^{i}du^{j}$ of $S.$ Let $\varphi (u^{1},u^{2})$ be a sufficient differentiable function on $M^{2}$. Then the second differential parameter of Beltrami with respect to the second fundamental form of $M^{2}$ is defined by \cite{K1}

\begin{equation}  \label{3}
\Delta ^{II}\varphi :=-\frac{1}{\sqrt{|b|}}(\sqrt{|b|}b^{ij}\varphi _{/i})_{/j},%
\end{equation}
where $\left( b^{ij}\right)$ denotes the inverse tensor of $\left(b_{ij}\right) $ and $b:=\det (b_{ij}).$ Applying (\ref{3}) for the position vector $\boldsymbol{x}$ of $M^{2}$, we find
\begin{equation}  \label{4}
\Delta ^{II}\boldsymbol{x} =-\frac{1}{2K}\nabla^{III}(K,\boldsymbol{n})-2\boldsymbol{n}%
\end{equation}
where $K, \boldsymbol{n}$ and $\nabla^{III}$ denote the curvature, the unit normal vector field and the first Beltrami-operator with respect to third fundamental form $III$. From (\ref{4}) we obtain the following results which were proved in \cite{S1}.
\begin{theorem}
\label{T1} A surface $M^{2}$ in $\mathbb{E}^{3}$ is of finite $II$-type 1 if and only if $M^{2}$ is part of a sphere.
\end{theorem}
\begin{theorem}
\label{T2} The Gauss map of a surface $M^{2}$ in $\mathbb{E}^{3}$ is of finite $II$-type 1 if and only if $M^{2}$ is part of a sphere.
\end{theorem}
Up to now, the only known surfaces of finite $II$-type in $\mathbb{E}^{3}$ are parts of spheres. So the following question seems to be interesting:

\begin{problem}
Other than the spheres, which surfaces in $\mathbb{E}^{3}$ are of finite $II$-type?
\end{problem}

%Another generalization of the above problem is to study surfaces in $\mathbb{E}^{3}$ of coordinate finite type, that is, their position vector $\boldsymbol{x}$ satisfying the relation
%\begin{equation}
%\Delta ^{II}\boldsymbol{x}=A\boldsymbol{x,}  \label{5}
%\end{equation}%
%where $A\in \mathbb{Re}^{3\times 3}$.
%
%From this point of view, we also pose the following problem
%
%\begin{problem}
%Classify all surfaces in $\mathbb{E}^{3}$ with the position vector $\boldsymbol{x}$ satisfying relation (\ref{5}).
%\end{problem}

This paper provides the first attempt at the study of finite type Gauss map of surfaces in $\mathbb{E}^{3}$ corresponding to the second fundamental form. In general when the Gauss map of $M^{2}$ is of finite type $k$, then there exist a monic polynomial, say $R(x)\neq 0,$ such that $R(\Delta ^{II})(\boldsymbol{n}-\boldsymbol{c})=\mathbf{0}.$ Suppose that $R(x)=x^{k}+\sigma_{1}x^{k-1}+...+\sigma _{k-1}x+\sigma _{k},$ then coefficients $\sigma _{i}$ are given by

\begin{eqnarray}
\sigma _{1} &=&-(\lambda _{1}+\lambda _{2}+...+\lambda _{k}),  \notag \\
\sigma _{2} &=&(\lambda _{1}\lambda _{2}+\lambda _{1}\lambda_{3}+...+\lambda
_{1}\lambda _{k}+\lambda _{2}\lambda _{3}+...+\lambda _{2}\lambda
_{k}+...+\lambda _{k-1}\lambda _{k}),  \notag \\
\sigma _{3} &=&-(\lambda _{1}\lambda _{2}\lambda _{3}+...+\lambda
_{k-2}\lambda _{k-1}\lambda _{k}),  \notag \\
&&.............................................  \notag \\
\sigma _{k} &=&(-1)^{k}\lambda _{1}\lambda _{2}...\lambda _{k}.  \notag
\end{eqnarray}

Therefore the position vector $\boldsymbol{n}$ satisfies the following equation, (see \cite{C3})

\begin{equation}  \label{2}
(\Delta^{II})^{k}\boldsymbol{n}+\sigma_{1}(\Delta^{II})^{k-1}\boldsymbol{n}+...+\sigma_{k}( \boldsymbol{n}-\boldsymbol{c})=\boldsymbol{0}.
\end{equation}

Our main result is the following

\begin{theorem} \label{T3}
All tubes in $\mathbb{E}^{3}$ are of infinite type corresponding to the second fundamental form.
\end{theorem}
Our discussion is local, which means that we show in fact that any open part of the Gauss map of a tube is of infinite Chen type.

\section{Tubes in $\mathbb{E}^{3}$}

Let $\ell: \boldsymbol{\varrho} =\boldsymbol{\varrho}(\emph{u})$, $\mathit{u}\epsilon (a,b)$ be a regular unit speed curve of finite length which is topologically imbedded in $\mathbb{E}^{3}$. The total space $N_{\boldsymbol{\varrho}}$ of the normal bundle of $\boldsymbol{\varrho}((a, b))$ in $\mathbb{E}^{3}$ is naturally diffeomorphic to the direct product $(a,b)\times \mathbb{E}^{2}$ via the translation along $\boldsymbol{\varrho}$ with respect to the induced normal connection. For a sufficiently small $r>0$ the tube of radius $r$ about the curve $\boldsymbol{\varrho}$ is the set:
\begin{equation}
T_{r}( \boldsymbol{\varrho})=\{exp_{\boldsymbol{\varrho}(u)}\boldsymbol{w}\mid\boldsymbol{w}\in N_{\boldsymbol{\varrho}} , \ \ \parallel \boldsymbol{w}%
\parallel =r,\ \ u\in(a,b)\}.  \notag
\end{equation}
Assume that ${\boldsymbol{t}, \boldsymbol{h}, \boldsymbol{b}}$ is the Frenet frame and $\kappa$ the curvature of the unit speed curve $\boldsymbol{\varrho} =
\boldsymbol{\varrho}(\emph{u})$. For a small real number $r$ satisfies $0 < r < min\frac{1}{|\kappa|}$, the tube $T_{r}( \boldsymbol{\varrho})$ is a smooth
surface in $\mathbb{E}^{3}$, \cite{R12}. Then a parametric representation of the tube $T_{r}( \boldsymbol{\varrho})$ is given by

\begin{equation}  \label{eq6}
\mathfrak{F}:\boldsymbol{x}(u,\varphi)= \boldsymbol{\varrho}+ r \cos\varphi\boldsymbol{h}+ r \sin\varphi\boldsymbol{b}.
\end{equation}

It is easily verified that the first and the second fundamental forms of $\mathfrak{F}$ are given by
\begin{align*}
I &= \big(\delta^{2}+r^{2}\tau^{2}\big)du^{2} + 2r^{2}\tau dud\varphi+
r^{2}d\varphi^{2}, \\
II &= \big(-\kappa\delta\cos\varphi+r\tau^{2}\big)du^{2} + 2r\tau dud\varphi
+rd\varphi^{2},
\end{align*}
where $\delta: = (1-r\kappa\cos\varphi)$ and $\tau$ is the torsion of the curve $\boldsymbol{\varrho}$. The Gauss curvature of $\mathfrak{F}$ is given by

\begin{equation}  \label{eq7}
K =-\frac{\kappa\cos\varphi}{r\delta}.
\end{equation}

Notice that $\kappa\neq0$ since the Gauss curvature vanishes. The Beltrami operator corresponding to the second fundamental form of $\mathfrak{F}$ can be
expressed as follows

\begin{eqnarray}  \label{eq8}
\Delta ^{II} &=&\frac{1}{\kappa\delta \cos \varphi }\Bigg[\frac{\partial^{2}}{\partial u^{2}}-2\tau \frac{\partial ^{2}}{\partial u\partial \varphi }%
+\Big(\tau ^{2}-\frac{\kappa\delta \cos \varphi}{r}\Big)\frac{\partial^{2}}{\partial \varphi^{2}}
\\
& &+\frac{(1-2\delta)\beta }{2\kappa\delta \cos \varphi }\frac{\partial }{\partial u}+\left(-\tau \acute{} +\frac{\tau \beta(2\delta-1) }{2\kappa\delta \cos \varphi }+\frac{\kappa(2\delta-1)\sin\varphi}{2r}\right) \frac{\partial }{\partial \varphi }\Bigg] ,
\notag
\end{eqnarray}
where $\beta: =\kappa \acute{}\cos \varphi +\kappa \tau \sin \varphi $ and $\acute{}:=\frac{d}{du}$.

The Gauss map $\boldsymbol{n}(u,\varphi)$ of $\mathfrak{F}$ has parametric representation
\begin{equation}  \label{eq81}
\boldsymbol{n}(u,\varphi)= -\cos \varphi \boldsymbol{h} - \sin \varphi \boldsymbol{b}.
\end{equation}

Applying (\ref{eq8}) for the position vector $\boldsymbol{n}$, we obtain
\begin{equation}  \label{eq82}
\Delta ^{II} \boldsymbol{n}=\frac{\beta}{2\kappa\delta^{2}\cos \varphi}\boldsymbol{t}+ \bigg(\frac{\sin^{2}\varphi}{2r\delta \cos\varphi}+\frac{\cos \varphi}{r\delta}-\frac{2\cos \varphi}{r}\bigg)\boldsymbol{h}+\frac{(1-4\delta)\sin\varphi}{2r\delta}\boldsymbol{b},
\end{equation}
which can be written
\begin{equation}  \label{eq83}
\Delta ^{II} \boldsymbol{n}=\frac{\beta}{2\kappa\delta^{2}\cos \varphi}\boldsymbol{t}+ \frac{1}{\kappa\delta \cos \varphi}\boldsymbol{P_{1}}(\cos \varphi, \sin\varphi),
\end{equation}
where $\boldsymbol{P_{1}}(\cos \varphi, \sin\varphi)$ is a vector with components polynomials of the functions $\cos \varphi$ and $\sin \varphi$ with coefficients functions of the variable $u$.

After long calculations, $(\Delta ^{II})^{2} \boldsymbol{n}$ can be computed as follows
\begin{equation}  \label{eq84}
(\Delta ^{II})^{2} \boldsymbol{n}=\frac{(3\delta-2)(12\delta-7)\beta^{3}}{4\kappa^{4}\delta^{5}\cos^{4} \varphi}\boldsymbol{t}+ \frac{1}{(\kappa\delta \cos \varphi)^{4}}\boldsymbol{P_{2}}(\cos \varphi, \sin\varphi),
\end{equation}
where $\boldsymbol{P_{2}}(\cos \varphi, \sin\varphi)$ is a vector with components polynomials of the functions $\cos \varphi$ and $\sin \varphi$ with coefficients functions of the variable $u$.

We need the following lemma which can be proved directly by using (\ref{eq8}).

\begin{lemma}
\label{L1} For any natural numbers $m$ and $n$, %and for any polynomial $h$ of the variable $\delta$ of degree $d$
we have
\begin{equation}
\Delta^{II} \Bigg(\frac{h(\delta)\beta^{m}}{\delta^{n}(\kappa\cos\varphi)^{n-1}}\Bigg)= -\frac{\widetilde{h}(\delta)\beta^{m+2}}{2\delta^{n+3}(\kappa\cos\varphi)^{n+2}}+\frac{1}{(\kappa\delta\cos\varphi)^{n+2}}Q(\cos\varphi,\sin\varphi),  \notag
\end{equation}
\end{lemma}
\noindent where $h(\delta)$ is a polynomial in $\delta$ of degree $d$, $Q$ is a polynomial in $\cos\varphi, \sin\varphi$ of degree $n$ + 3 with functions in $u$ as coefficients, deg($\widetilde{h}(\delta)$) = $d$ +2 and
\begin{equation}
\widetilde{h}(\delta) = ((2n-1)\delta -n)(4(n+1)\delta -(2n+3))h(\delta).
\end{equation}

So, it is easily verified that
\begin{equation}  \label{eq85}
(\Delta ^{II})^{3} \boldsymbol{n}=\frac{(3\delta-2)(12\delta-7)(9\delta-7)(24\delta-13)\beta^{5}}{8\kappa^{7}\delta^{8}\cos^{7} \varphi}\boldsymbol{t}+ \frac{1}{(\kappa\delta \cos \varphi)^{7}}\boldsymbol{P_{3}}(\cos \varphi, \sin\varphi),
\end{equation}
where $\boldsymbol{P_{3}}(\cos \varphi, \sin\varphi)$ is a vector with components polynomials of the functions $\cos \varphi$ and $\sin \varphi$ with coefficients functions of the variable $u$.

Moreover, by using Lemma (\ref{L1}), one can find
\begin{equation}  \label{eq86}
(\Delta ^{II})^{k} \boldsymbol{n}=\frac{h_{k}(\delta)\beta^{2k-1}}{2^{k}\delta^{3k-1}(\kappa\cos \varphi)^{3k-2}}\boldsymbol{t}+ \frac{1}{(\kappa\delta \cos \varphi)^{3k-2}}\boldsymbol{P_{k}}(\cos \varphi, \sin\varphi),
\end{equation}
and
\begin{equation}  \label{eq87}
(\Delta ^{II})^{k+1} \boldsymbol{n}=\frac{h_{k+1}(\delta)\beta^{2k+1}}{2^{k+1}\delta^{3k+2}(\kappa\cos \varphi)^{3k+1}}\boldsymbol{t}+ \frac{1}{(\kappa\delta \cos \varphi)^{3k+1}}\boldsymbol{P_{k+1}}(\cos \varphi, \sin\varphi),
\end{equation}
where $\boldsymbol{P_{k}}(\cos \varphi, \sin\varphi)$, $\boldsymbol{P_{k+1}}(\cos \varphi, \sin\varphi)$ are vectors with components polynomials of the functions $\cos \varphi, \sin \varphi$ with coefficients functions of the variable $u$, and
\begin{equation}
h_{\lambda}(\delta) = \overset{\lambda-1}{\underset{j=1}{\prod }}12j\delta(4j-5) \neq 0,  \notag
\end{equation}
for any positive integer $\lambda$.

Now, we suppose that the tube $\mathfrak{F}$ is of finite type, thus from (\ref{eq83}-\ref{eq87}), and (\ref{2})we find that
\begin{eqnarray}  \label{eq88}
\frac{h_{k+1}(\delta)\beta^{2k+1}}{2^{k+1}\delta^{3k+2}(\kappa\cos \varphi)^{3k+1}}\boldsymbol{t}+ \frac{1}{(\kappa\delta \cos \varphi)^{3k+1}}\boldsymbol{P_{k+1}} &&  \notag \\
+ \frac{\sigma_{1} h_{k}(\delta)\beta^{2k-1}}{2^{k}\delta^{3k-1}(\kappa\cos \varphi)^{3k-2}}\boldsymbol{t}+ \frac{\sigma_{1}}{(\kappa\delta \cos \varphi)^{3k-2}}\boldsymbol{P_{k}}&&  \notag \\
+ ...+ \frac{\sigma_{k-1}(3\delta-2)(12\delta-7)\beta^{3}}{4\kappa^{4}\delta^{5}\cos^{4} \varphi}\boldsymbol{t}+ \frac{\sigma_{k-1}}{(\kappa\delta \cos \varphi)^{4}}\boldsymbol{P_{2}} &&  \notag \\
+\frac{\sigma_{k}\beta}{2\kappa\delta^{2}\cos \varphi}\boldsymbol{t}+ \frac{\sigma_{k}}{\kappa\delta \cos \varphi}\boldsymbol{P_{1}} &=& \boldsymbol{0}. \notag
\end{eqnarray}

The above equation can be simply written as follows
\begin{equation}  \label{eq89}
\frac{h_{k+1}(\delta)\beta^{2k+1}}{\delta} \boldsymbol{t}= Q_{1}(\cos \varphi,\sin\varphi)\boldsymbol{t}+Q_{2}(\cos \varphi,\sin\varphi)\boldsymbol{h}+ Q_{3}(\cos \varphi,\sin\varphi)\boldsymbol{b},
\end{equation}
where $Q_{i}(\cos \varphi,\sin\varphi)$ are polynomials in $\cos \varphi$ and $\sin\varphi$.
%Accordingly, we get

We distinguish two cases

Case I. $\beta \neq 0$. From (\ref{eq89}), we have
\begin{equation}  \label{eq90}
\frac{h_{k+1}(\delta)\beta^{2k+1}}{\delta} = Q_{1}(\cos \varphi,\sin\varphi).
\end{equation}
This is impossible for any $k \geq 1$ since $h_{k+1}(\delta)\neq 0$.

Case II. $\beta \equiv 0$. Then $\kappa' = 0$ and $\kappa\tau = 0$. Thus $\kappa$ = const. $\neq 0$ and $\tau = 0$, therefore the curve $\boldsymbol{\varrho}$ is a plane circle and so, $\mathfrak{F}$ is an anchor ring.

%A parametric representation of an anchor ring can be expressed as follows
%\begin{equation}  \label{eq9}
%\mathfrak{F}:\boldsymbol{x}(u,\varphi)= \{(a+ r\cos u)\cos\varphi,(a+ r\cos u)\sin\varphi, r\sin u \},
%\end{equation}
%\begin{center}
%$a > r,\ \ a \epsilon \mathbb{R}.$
%\end{center}

In this case, the first fundamental form becomes

\begin{equation*}
I=\delta^{2}du^{2}+r^{2}d\varphi ^{2},
\end{equation*}%
while the second is

\begin{equation*}
II= -\kappa\delta\cos\varphi du^{2} +rd\varphi^{2}.
\end{equation*}

Hence, equation (\ref{eq8}) reduces to
\begin{equation}  \label{eq10}
\Delta^{II}= \frac{1}{\kappa \delta \cos \varphi}\bigg(\frac{\partial^{2}}{\partial u^{2}}-\frac{\kappa \delta \cos \varphi}{r}\frac{\partial^{2}}{\partial\varphi^{2}} +\frac{\kappa(2\delta-1)\sin \varphi}{2r}\frac{\partial}{\partial \varphi}\bigg).
\end{equation}

Applying (\ref{eq10}) for the position vector $\boldsymbol{n}$, one finds
\begin{equation}  \label{eq91}
\Delta ^{II} \boldsymbol{n}=\bigg(\frac{\sin^{2}\varphi}{2r\delta \cos\varphi}+\frac{\cos \varphi}{r\delta}-\frac{2\cos \varphi}{r}\bigg)\boldsymbol{h}+\frac{(1-4\delta)\sin\varphi}{2r\delta}\boldsymbol{b},
\end{equation}
which can be written as follows
\begin{equation}  \label{eq11}
\Delta ^{II} \boldsymbol{n}=\frac{\sin^{2}\varphi}{2r\delta \cos\varphi}\boldsymbol{h}+\frac{1}{\delta}\boldsymbol{F_{1}}(\cos \varphi, \sin\varphi).
\end{equation}

Consequently, we get
\begin{equation}  \label{eq12}
(\Delta ^{II})^{2} \boldsymbol{n}=-\frac{3\sin^{4}\varphi}{4r^{2}\delta^{3} \cos^{3}\varphi}\boldsymbol{h}+\frac{1}{\delta^{3} \cos^{2}\varphi}\boldsymbol{F_{2}}(\cos \varphi, \sin\varphi),
\end{equation}
where $\boldsymbol{F_{1}}(\cos \varphi, \sin\varphi)$, $\boldsymbol{F_{2}}(\cos \varphi, \sin\varphi)$ are vectors with components polynomials of the functions $\cos \varphi, \sin \varphi$ with coefficients functions of the variable $u$.

One can easily prove that
\begin{equation}  \label{eq13}
\Delta ^{II} \frac{\sin^{m}\varphi}{(\delta \cos\varphi)^{n}} =\frac{3\sin^{m+2}\varphi}{2r(\delta \cos\varphi)^{n+2}}+\frac{1}{(\delta\cos\varphi)^{n+1}}Q(\cos \varphi, \sin\varphi),
\end{equation}
where $Q(\cos \varphi, \sin\varphi)$ is a polynomial in $\cos \varphi, \sin \varphi$ with coefficients functions of the variable $u$. Therefore, we find that

\begin{equation}  \label{eq14}
(\Delta ^{II})^{k} \boldsymbol{n}=(-1)^{k-1}\frac{3^{k-1}\sin^{2k}\varphi}{(2r)^{k}(\delta \cos \varphi)^{2k-1}}\boldsymbol{h}+\frac{1}{\delta^{2k-1} \cos ^{2k-2}\varphi}\boldsymbol{F_{k}}(\cos \varphi, \sin\varphi),
\end{equation}
and
\begin{equation}  \label{eq15}
(\Delta ^{II})^{k+1} \boldsymbol{n}=(-1)^{k}\frac{3^{k}\sin^{2k+2}\varphi}{(2r)^{k+1}(\delta \cos \varphi)^{2k+1}}\boldsymbol{h}+\frac{1}{\delta^{2k+1} \cos ^{2k}\varphi}\boldsymbol{F_{k+1}}(\cos \varphi, \sin\varphi),
\end{equation}
where $\boldsymbol{F_{k}}(\cos \varphi, \sin\varphi)$, $\boldsymbol{F_{k+1}}(\cos \varphi, \sin\varphi)$ are vectors with components polynomials of the functions $\cos \varphi, \sin \varphi$ with coefficients functions of the variable $u$.

On account of (\ref{eq11}), (\ref{eq12}), (\ref{eq14}), (\ref{eq15}) and (\ref{2}), we conclude that
\begin{eqnarray}  \label{eq16}
(-1)^{k}\frac{3^{k}\sin^{2k+2}\varphi}{(2r)^{k+1}(\delta \cos \varphi)^{2k+1}}\boldsymbol{h}+\frac{1}{\delta^{2k+1} \cos ^{2k}\varphi}\boldsymbol{F_{k+1}}(\cos \varphi, \sin\varphi) &&  \\
+ (-1)^{k-1}\sigma_{1}\frac{3^{k-1}\sin^{2k}\varphi}{(2r)^{k}(\delta \cos \varphi)^{2k-1}}\boldsymbol{h}+\sigma_{1}\frac{1}{\delta^{2k-1} \cos ^{2k-2}\varphi}\boldsymbol{F_{k}}(\cos \varphi, \sin\varphi) &&  \notag \\
+ ...-\sigma_{k-1}\frac{3\sin^{4}\varphi}{4r^{2}\delta^{3} \cos^{3}\varphi}\boldsymbol{h}+\sigma_{k-1}\frac{1}{\delta^{3} \cos^{2}\varphi}\boldsymbol{F_{2}}(\cos \varphi, \sin\varphi) &&  \notag \\
+ \sigma_{k}\frac{\sin^{2}\varphi}{2r\delta \cos\varphi}\boldsymbol{h}+\sigma_{k}\frac{1}{\delta}\boldsymbol{F_{1}}(\cos \varphi, \sin\varphi) &=& \boldsymbol{0}. \notag
\end{eqnarray}

For simplicity, equation (\ref{eq16}) can be written as follows
\begin{eqnarray}  \label{eq17}
\frac{3^{k}\sin^{2k+2}\varphi}{(2r)^{k+1}\cos \varphi}\boldsymbol{h}=Q_{1}(\cos \varphi, \sin\varphi)\boldsymbol{h}+ Q_{2}(\cos \varphi, \sin\varphi)\boldsymbol{b},
\end{eqnarray}
where $Q_{i}(\cos \varphi, \sin\varphi)$, $i =1,2$ are polynomials in $\cos \varphi, \sin \varphi$ with coefficients functions of the variable $u$. Therefore, we find that
\begin{eqnarray}  \label{eq18}
\frac{3^{k}\sin^{2k+2}\varphi}{(2r)^{k+1}\cos \varphi}= Q_{1}(\cos \varphi, \sin\varphi).
\end{eqnarray}

This is impossible, since the right side hand $Q_{1}$ is a polynomial in $\cos \varphi, \sin \varphi$ while the left side hand is not. Thus, our theorem is proved.
%we have just proved the following

%\begin{corollary}
%\label{C1.1} The Gauss map of an anchor ring in the Euclidean 3-space is of infinite $II$-type.
%\end{corollary}
%
%\begin{theorem}
%\label{T3} The Gauss map of a tube in $\mathbb{E}^{3}$ is of infinite type corresponding to the second fundamental form.
%\end{theorem}

%It can be seen that $d_{\lambda }\neq 0$, for each natural number $\lambda $. Moreover, we have

%\begin{acknowledgement}
%The authors would like to express their thanks to the referee for his useful remarks.
%\end{acknowledgement}

%\bibliography{...}

\end{document}